\documentclass[pdflatex,sn-mathphys-num]{sn-jnl}
\usepackage{graphicx}
\usepackage{amsmath,amssymb,amsfonts}%
\usepackage[title]{appendix}%

\usepackage{mathtools}

\usepackage{bm}
\usepackage{cleveref}
\usepackage{enumitem}
\usepackage{algorithm}
\usepackage{siunitx}
\usepackage{subfig}

\newcommand{\bu}{\bm{u}}
\newcommand{\bq}{\bm{q}}
\newcommand{\br}{\bm{r}}
\newcommand{\bv}{\bm{v}}
\newcommand{\bx}{\bm{x}}
\newcommand{\bn}{\bm{n}}
\newcommand{\myd}{\mathrm{d}}
\newcommand{\Mat}[1]{M_{\text{\tiny{#1}}}}

\theoremstyle{thmstyleone}%
\theoremstyle{thmstyletwo}%
\newtheorem{remark}{Remark}%

\theoremstyle{thmstylethree}%

\numberwithin{equation}{section}
\begin{document}

\title[Solving Biot with OPM Flow and TPSA]{Solving Biot poroelasticity by coupling OPM Flow with the two-point stress approximation finite volume method}
\author*[1]{\fnm{Wietse Marijn} \sur{Boon}}\email{wibo@norceresearch.no}
\author[1]{\fnm{Sarah} \sur{Gasda}} 
\author[1]{\fnm{Tor Harald} \sur{Sandve}} 
\author[1]{\fnm{Svenn} \sur{Tveit}}

\affil[1]{\orgdiv{Division of Energy and Technology}, \orgname{NORCE Norwegian Research Centre}, \orgaddress{\street{Nygårdsgaten 112}, \city{Bergen}, \postcode{5008}, \country{Norway}}}

\abstract{
    Finite volume methods are prevalent in reservoir simulation due to their mass conservation properties and their ability to handle complex grids. 
    However, a simple and consistent finite volume method for elasticity was unavailable until the recently developed two-point stress approximation finite volume method (TPSA). In this work, we show how to couple TPSA to an established flow simulator, using OPM Flow as our primary example. Due to this choice of numerical methods, the coupling is naturally handled at the cell centers, without requiring interpolation operators. We propose a fixed stress coupling scheme and reuse algebraic multi-grid preconditioners, which are known to be effective for two-point flux finite volume methods. Numerical examples illustrate the flexibility of the approach and we showcase how the introduction of solid mechanics impacts the behavior of compartmentalized flow systems.
}

\maketitle

\section{Introduction}

The interaction between fluid flow and solid mechanics plays an important role in reservoir engineering applications such as CO$_2$ storage, in which the geological site may deform due to artificially induced pressure changes. Simulating such systems efficiently and accurately remains a computational challenge, especially on the regional scale. 

The choice of numerical methods is one of the partial causes for this difficulty. By convention, finite volume methods are predominantly chosen to simulate fluid flow, whereas solid mechanics is typically modeled using the finite element method. However, finite element methods are more sensitive to the aspect ratios of the grid cells and often struggle to handle the complex corner-point grids commonly used in reservoir models. Moreover, finite element and finite volume methods place the degrees of freedom at different locations in the grid, which means that interpolation operators need to be employed to communicate between the two.

Numerous discretization techniques have been proposed in the literature to solve the coupled poroelasticity equations. For example, finite element methods based on the primal formulation of elasticity in primal form were investigated in \cite{phillips2007coupling,rodrigo2016stability,riviere2017error}, and these were coupled to finite volume methods for flow in \cite{pettersen2009improved,settgast2024geos}. Mixed finite elements for the elasticity system were investigated, e.g., in \cite{lee2016robust,caucao2022multipoint} and alternatives have been proposed based on virtual element methods \cite{kumar2024numerical,botti2025fully}, hybridizable discontinuous Galerkin methods \cite{fu2019high,cesmelioglu2024hybridizable} and polygonal discontinuous Galerkin methods \cite{zhao2023locking}, as well as rotation-based formulations of elasticity \cite{anaya2020rotation,boon2023mixed}. 

For reservoir simulators based on finite volume methods, it is attractive to employ a finite volume discretization for the solid mechanics as well. However, their extension to the elasticity equations is non-trivial. Finite volume methods place the degrees of freedom in the cell centers and approximate derivatives using finite difference stencils \cite{leveque2002finite}. In the case of elasticity, we can evaluate the normal gradient of the displacement using a two-point finite difference, but to compute the normal component of the \emph{symmetric} gradient, one typically requires a larger stencil. Such stencils are employed in the staggered grid approach of the MAC scheme \cite{harlow1965numerical} or the Multi-Point Stress Approximation method \cite{nordbotten2016stable,keilegavlen2017finite}.
The recently introduced Two-Point Stress Approximation (TPSA) method circumvents this issue by introducing two additional variables to the problem, namely a solid pressure and a rotation, to form a consistent finite volume method with 7 degrees of freedom per cell, and a minimal stencil \cite{nordbotten2025two}.

In this work, we show how the TPSA method can be coupled to a reservoir simulator such as OPM Flow \cite{rasmussen2021open} using a fixed stress splitting scheme \cite{mikelic2013convergence,both2017robust}.
We highlight the benefits of employing two finite volume methods on the same grid. In particular, the information passed between the flow and elasticity solvers consists of only one value per cell, without the need for interpolation operators. Moreover, optimized solvers designed for the two-point flux (TPFA) finite volume method can be reused to efficiently solve the TPSA system.

We use the coupled model to investigate the secondary effects that poroelasticity can have on reservoir models. In such models, impermeable zones or regions with low pore volumes are often disregarded since they prevent fluid flow. This may cause the model to decouple into separate flow systems that evolve independently through time. However, impermeable zones still play a connecting role between the compartments through the solid mechanics equations. We illustrate this phenomenon in an isolated setting in \Cref{sub:a_sealing_barrier}. In a simplified example, we show that poroelasticity can affect the fluid pressure through sealing barriers. This secondary impact on the fluid pressure evolution indicates the significance of incorporating solid mechanics in reservoir simulations.

The remainder of this article is organized as follows. First, \Cref{sub:the_poroelasticity_model} introduces the linear Biot model for poroelasticity. \Cref{sec:TPSA} provides a brief description of the TPSA method for linearized elasticity. The fixed stress splitting scheme is described in \Cref{sec:fixed stress split} and we highlight the details regarding implementation in \Cref{sec:implementation}. \Cref{sec:numerical_results} presents the numerical results and we present the conclusions in \Cref{sec:conclusions}.

\subsection{The poroelasticity model}
\label{sub:the_poroelasticity_model}

Let us briefly recall the linearized Biot model in the case of single-phase flow. 
First, we assume that the solid mechanics are governed by Hooke's law and conservation of linear momentum:
\begin{subequations} \label{eqs: Biot original}
\begin{align}
    \sigma &= 2 \mu \varepsilon(\bu) + \lambda (\nabla \cdot \bu) I, \label{eq:Cauchy}\\
	- \nabla \cdot(\sigma - \alpha (p_f - p_0) I) &= \bm{f_u}. \label{eq: momentum balance orig}
\end{align}
Here, $\sigma$ represents the Cauchy stress tensor and $\bu$ the displacement of the poroelastic medium. $\varepsilon$ denotes the symmetric gradient, $\mu$ and $\lambda$ are the Lamé parameters, and $\alpha$ is the Biot-Willis constant. $p_f - p_0$ is the deviation of the fluid pressure with respect to a given reference pressure $p_0$. This term acts as an isotropic stress in the momentum balance equation. $\bm{f_u}$ contains body forces.

Second, the fluid dynamics are governed by Darcy's law and mass conservation:
\begin{align}
    \bq &= - \frac{K}{\mu_w} (\nabla p_f - \rho \bm{g}), \\
    \partial_t (c_0 p_f + \alpha \nabla \cdot \bu) + \nabla \cdot \bq &= f_p \label{eq: mass balance orig}
\end{align}
\end{subequations}
in which $K$ is the permeability of the solid, $\mu_w$ the viscosity and $\rho$ the density of the fluid, $\bm{g}$ the gravity force, and $c_0$ the storativity. $\partial_t$ denotes the partial derivative with respect to time.
We note that the term $\partial_t (\alpha \nabla \cdot \bu)$ models how volumetric changes in the solid affect the fluid pressure. 

\section{The Two-Point Stress Approximation method for linearized elasticity}
\label{sec:TPSA}

In this section, we provide a brief summary of the two-point stress approximation method (TPSA), proposed in \cite{nordbotten2025two}.
The method is based on a reformulation of the elasticity equations in terms of three primary variables. 
In addition to the displacement $\bu$, it introduces a rotation variable $\br$ and a solid pressure $p_s$:
\begin{align} \label{eq: r and p_s}
    \br &= \mu \nabla \cdot (S^* \bu), &
    p_s &= \lambda \nabla \cdot \bu.
\end{align}
Here, $S^*$ is the ``skew'' operator that maps 3-vectors to anti-symmetric $3 \times 3$ matrices:
\begin{align}
    S^* \bu \coloneqq 
    \begin{bmatrix}
        0 & -u_3 & u_2 \\
        u_3 & 0 & -u_1 \\
        - u_2 & u_1 & 0
    \end{bmatrix}.
\end{align}
We remark that $S^*$ is related to the cross product by the identity $(S^* \bu)\bv = \bu \times \bv$.

The Cauchy stress \eqref{eq:Cauchy} can now be rewritten as $\sigma = 2\mu (\nabla \bu) + S^* \br + p_s I$ by using the identity $S^* \nabla \cdot (S^* \bu) = (\nabla \bu)^T - \nabla \bu$.
We define two additional dual variables, namely $\tau = S^* \bu$ and $\bv = \bu$. 
While these variables are superfluous in the continuous setting, they are handled differently by the discretization.
To summarize, the dual variables are given by 
\begin{subequations} \label{eq:TPSA_system_cont}
\begin{align} \label{eq: dual var map}
    \begin{bmatrix}
        \sigma \\ \tau \\ \bv
    \end{bmatrix}
    =
    \begin{bmatrix}
        2 \mu \nabla & S^* & I \\
        S^* \\
        I 
    \end{bmatrix}
    \begin{bmatrix}
        \bu \\ \br \\ p_s
    \end{bmatrix},
\end{align}
and the system of equations governing linearized elasticity become
\begin{align} \label{eq: conservation 3var}
    -\nabla \cdot 
    \begin{bmatrix}
        \sigma \\ \tau \\ \bm{v}
    \end{bmatrix}
    + 
    \begin{bmatrix}
        0 \\ 
        & \mu^{-1} \\ 
        & & \lambda^{-1} \\
    \end{bmatrix}
    \begin{bmatrix}
        \bu \\ \br \\ p_s
    \end{bmatrix}
    =
    \begin{bmatrix}
        \bm{f_u} \\ 0 \\ 0
    \end{bmatrix}
\end{align}
\end{subequations}

\subsection{Discretization}
\label{sub:TPSA_discretization}

The TPSA method discretizes \eqref{eq:TPSA_system_cont} by placing the primary unknowns $(\bu, \br, p_s)$ in the cell centers and the dual variables $(\sigma, \tau, \bm{v})$ on the faces of the mesh. Similar to the two-point flux (TPFA) finite volume method for flow, we only consider the normal components of the dual variables on the faces. For example, $\sigma_k$ is a 3-vector that denotes the traction on face $\varsigma_k$.

The discretization of equation \eqref{eq: dual var map} involves taking weighted averages of the primal variables onto the faces. These are defined as follows. Let an interior face $\varsigma_k$ be between cells $\omega_i$ and $\omega_j$ such that its unit normal vector $\bn_k$ coincides with $\bn_i$, which is outward with respect to $\omega_i$. Let $\delta_{ik} \coloneqq \bn_i \cdot (\bx_k - \bx_i)$ be the normal distance between face center $\bx_k$ and cell center $\bx_i$. Then the weights and averaging operators are given by
\begin{align}
    w_i &\coloneqq \frac{\delta_{ik}}{\mu_i}, &
    \widetilde \Xi_k \bu &\coloneqq \frac{w_i \bu_i + w_j \bu_j}{w_i + w_j}, &
    \Xi_k \bu &\coloneqq \frac{w_j \bu_i + w_i \bu_j}{w_i + w_j},
\end{align}
with $\mu_i$ the elastic modulus in cell $i$. Let $\nabla_k \bu \coloneqq \frac{\bu_j - \bu_i}{\delta_k}$ be the conventional two-point approximation of the normal derivative, with $\delta_k \coloneqq \delta_{ik} + \delta_{jk}$. Let $\bar \mu_k \coloneqq \widetilde \Xi \mu$ be the effective elasticity modulus.  The TPSA discretization of \eqref{eq: dual var map} on $\varsigma_k$ is then given by
\begin{align} \label{eq: dual var map discrete}
    \begin{bmatrix}
        \sigma_k \\ \tau_k \\ \bv_k
    \end{bmatrix}
    =
    |\varsigma_k|
    \begin{bmatrix}
        2 \bar \mu_k \nabla_k & -(S^*\bn_k) \widetilde \Xi_k & \bn_k \widetilde \Xi_k \\
        -(S^*\bn_k) \Xi_k \\
        \bn_k \Xi_k & & \delta_k^\mu \nabla_k
    \end{bmatrix}
    \begin{bmatrix}
        \bu \\ \br \\ p_s
    \end{bmatrix}
\end{align}

The term in the $(3,3)$-block of this operator is a parameter-free stabilization term. The weight is defined as $\delta_k^\mu \coloneqq \frac12 w_iw_j \bar \mu_k \sim h^2 \mu^{-1}$, i.e. it scales quadratically with the mesh size and inversely with $\mu$.

The discretization of \eqref{eq: conservation 3var}, on the other hand, follows by integrating over the cells and applying the divergence theorem.
In particular, we use the definition of the dual variables \eqref{eq: dual var map discrete} to derive the following discretization of \eqref{eq: conservation 3var} for a cell $\omega_i \in \Omega_h$:
\begin{align}
    \Mat{TPSA}^i \begin{bmatrix}
        \bu \\ \br \\ p_s
    \end{bmatrix}
    \coloneqq 
    - \sum_{\varsigma_k \subseteq \partial \omega_i} \epsilon_{ik}
    \begin{bmatrix}
        \sigma_k \\ \tau_k \\ \bm{v}_k
    \end{bmatrix}
    + 
    |\omega_i|
    \begin{bmatrix}
        0 \\ \mu_i^{-1} \br_i \\ \lambda_i^{-1} p_{s, i}
    \end{bmatrix},
\end{align}
in which $\epsilon_{ik} = \bn_i \cdot \bn_k = \pm 1$.
Finally, we collect the equations for all cells $\omega_i \in \Omega_h$ to obtain a linear system in terms of the primary unknowns, with 7 degrees of freedom per cell:
\begin{align} \label{eq:TPSA_system}
    \Mat{TPSA} \begin{bmatrix}
        \bu \\ \br \\ p_s
    \end{bmatrix}
    =
    \begin{bmatrix}
        |\omega| \bm{f_u} \\ 0 \\ 0
    \end{bmatrix}.
\end{align}

Boundary conditions are implemented by associating appropriate weights to the ``outside'' of boundary faces. To illustrate, let $\varsigma_k$ be a boundary face that borders cell $\omega_i$. To enforce homogeneous conditions on the boundary, we define $\bu_j \coloneqq 0$. A zero displacement, or fixed, boundary condition is realized by setting $\delta_{jk} = 0$ so that $\Xi_k \bu = 0$ and $\widetilde \Xi_k \bu = \bu_i$. Alternatively, a positive, bounded $\delta_{jk}$ and $\mu_j$ leads to Robin boundary conditions, which simulate springs with a spring constant $w_j^{-1}$. The implementation of zero traction, or free, boundaries follows by considering the limit of $\delta_{jk} \to \infty$. For more details, we refer to \cite{nordbotten2025two}.

\subsection{Relation to TPFA}

To complete the exposition, we briefly recall the two-point flux approximation finite volume method for Darcy flow, using the same notation.
On a mesh face $\varsigma_k = \partial \omega_i \cap \partial \omega_j$, we define the effective conductivity by the weighted harmonic average $\bar \kappa_k \coloneqq \frac{\delta_k}{\delta_{ik} \mu_{w, i} K_i^{-1} + \delta_{jk} \mu_{w, j} K_j^{-1}}$, similar to $\bar{\mu}_k$ in \eqref{eq: dual var map discrete}.
The normal flux is then approximated by
\begin{align}
    \bq_k = - |\varsigma_k| \left(\bar \kappa_k \nabla_k p_f - \rho \bn_k \cdot \bm{g} \right)
\end{align}
The accumulation of mass in cell $\omega_i$ is computed as the sum of fluxes:
\begin{align}
    \Mat{TPFA}^i (p_f) &\coloneqq \sum_{\varsigma_k \subseteq \partial \omega_i} \epsilon_{ik} \bq_k.
\end{align}
Collecting $\Mat{TPFA}^i$ for all cells, we form the linear system $\Mat{TPFA} p_f = |\omega| f_p$, which has one degree of freedom per cell. The TPFA finite volume method remains the industry standard because it leads to symmetric positive definite systems for which a range of efficient numerical solvers are available. The method is, however, not consistent in general and is only guaranteed to converge to the correct solution on so-called $K$-orthogonal grids \cite{aavatsmark2002introduction,eymard2000finite}. A similar consistency limitation holds for TPSA. Nevertheless, a major advantage is that TPFA and TPSA are both stable on the same, large class of grids.

\section{A splitting scheme using cell-centered variables}
\label{sec:fixed stress split}

In this section, we propose the splitting scheme to solve the poroelasticity problem using TPFA for the flow equations and TPSA for the mechanics. If we were to discretize the momentum balance \eqref{eq: momentum balance orig} directly, then we need to evaluate the term $\nabla \cdot (\alpha p_fI) = \alpha \nabla p_f$ at the cell centers. However, this is not immediately available because $p_f$ is a cell-centered variable and its gradient is more naturally evaluated on faces.

Similar to the reformulation of elasticity in \Cref{sec:TPSA}, we remedy this issue by reformulating the problem. Let us introduce the deviation of the fluid pressure as $\Delta p_f$ and the effective pressure $\hat p$ as
\begin{align} \label{eq: def eff pressure}
    \Delta p_f &\coloneqq p_f - p_0, &
    \hat p &\coloneqq \lambda \nabla \cdot \bu - \alpha \Delta p_f.
\end{align}
The Biot equations \eqref{eqs: Biot original} are then rewritten as:
\begin{subequations}
\begin{align}
    - \nabla \cdot(2 \mu \varepsilon(\bu) + \hat p I) &= \bm{f_u}, \label{eq:momentum balance}\\
    - \nabla \cdot \bu + \frac1{\lambda} \hat p + \frac{\alpha}{\lambda} \Delta p_f &= 0, \label{eq:def bar p}\\ 
    \left(c_0 + \frac{\alpha^2}{\lambda} \right) \partial_t p_f + \frac{\alpha}{\lambda} \partial_t \hat p + \nabla \cdot \bq &= f_p, \label{eq: mass balance}\\
    \bq + \frac{K}{\mu_w} \nabla p_f &= K \rho \bm{g}.
\end{align}
\end{subequations}

\begin{remark} \label{rem:compressibility}
    The reformulation introduces a term $\frac{\alpha^2}{\lambda}\partial_t p_f$ in the mass balance equation \eqref{eq: mass balance}. This term models the rock compressibility under the assumption of fixed effective pressure. In particular, if $\partial_t \hat p = 0$, then \eqref{eq: def eff pressure} implies $\partial_t \nabla \cdot \bu = \frac{\alpha}{\lambda}\partial_t p_f$. Substitution in \eqref{eq: mass balance orig} leads to exactly this term. A similar compressibility term appears in the fixed stress scheme of \cite{mikelic2013convergence,both2017robust}.
\end{remark}

\begin{remark}[Gravity]
    We assume that the poroelastic medium is in equilibrium at the beginning of the simulation, $\bu$ represents the displacement with respect to that reference configuration, and $p_0$ is the hydrostatic pressure distribution. $\Delta p_f = 0$ should then imply $\bu = 0$ and $\hat p = 0$, which is only the case if $\bm{f_u} = 0$. From this, we conclude that $\bm{f_u}$ corresponds to the body forces that are \emph{additional} to the gravity forces at $t = 0$, if present.
\end{remark}

Our splitting scheme is based on iteratively solving for $\hat p$ and $p_f$. In particular, if $\partial_t \hat p$ is known in \eqref{eq: mass balance}, then we may move it to the right-hand side and it will act as a mass source in the flow equations. Similarly, if $\Delta p_f$ is known in \eqref{eq:def bar p} and we move the related term to the right-hand side, then we recognize \eqref{eq:def bar p} as the third row of \eqref{eq: conservation 3var}. This leads us to the fixed stress scheme described in \Cref{alg: space-time}.
\begin{algorithm}[H]
    \caption{Fixed stress scheme}
    \label{alg: space-time}
\begin{enumerate}[leftmargin=*,label=\arabic*.] 
    \item Initialize $\hat p^0$ and set the iteration index $n = 1$.
    \item Use TPFA to solve the flow problem: find $p_f^n$ such that
    \begin{align} \label{eq: OPM solve}
        (c_0 + \frac{\alpha^2}{\lambda}) \partial_t p_f^n + \Mat{TPFA} p_f^n &= f_p - \frac{\alpha}{\lambda} \partial_t \hat p^{n - 1},
    \end{align}
    \item Use TPSA to solve the elasticity problem: find $[\bu^n, \br^n, \hat p^n]$ such that
        \begin{align} \label{eq: TPSA solve}
            \Mat{TPSA} \begin{bmatrix}
                \bu^n \\ \br^n \\ \hat p^n
            \end{bmatrix} =
            \begin{bmatrix}
                \bm{f_u} \\ 0 \\ - \frac{\alpha}{\lambda} \Delta p_f^n
            \end{bmatrix}
        \end{align}
    \item Increment $n$ and repeat the previous two steps until convergence.
\end{enumerate}
\end{algorithm}

A key advantage of this scheme is the minimal amount of information that needs to be passed between the flow and mechanics solvers. In particular, because both $p_f$ and $\hat p$ are scalar-valued, cell-wise variables, the solvers exchange a single value per grid cell at each time step and iteration.

We moreover emphasize that there is no need to interpolate, because the pressure variables $p_f$ and $\hat p$ are defined in the cell centers, where the right-hand side terms are evaluated. This is a simple consequence of coupling two finite volume methods, which significantly simplifies the implementation.

\subsection{Variations}
\label{sub:iteration_schemes}

In this subsection, we discuss variants of the splitting scheme illustrated in \Cref{alg: space-time}. 
We first recognize \Cref{alg: space-time} as a fixed point iterative scheme. 

Let $\psi = - \frac{\alpha}{\lambda} \partial_t \hat p$. Let $\mathcal{F}(\psi)$ be the operator that i) solves \eqref{eq: OPM solve} with $\psi$ on the right-hand side, ii) uses the computed fluid pressure to solve \eqref{eq: TPSA solve}, and iii) computes the new source term based on the computed total pressure.
\Cref{alg: space-time} then iteratively solves the fixed point problem $\mathcal{F}(\psi) = \psi$ by setting $\psi^{n + 1} = \mathcal{F}(\psi^n)$. 

The convergence of this fixed point scheme can be improved by employing Anderson acceleration \cite{walker2011anderson}. Instead of using the output $\mathcal{F}(\psi^n)$ as the new source term, Anderson acceleration involves taking a weighted arithmetic average of the previous $m = \min\{m_0, n\}$ iterates, for a pre-defined $m_0 \in \mathbb{N}$. The weights are chosen to minimize the residual based on the previous $m$ residuals. In particular, we solve the following least-squares problem
\begin{align} \label{eq: Anderson}
    &\min_{\bm{\beta} \in \mathbb{R}^m} \left\| \sum_{i = 0}^{m - 1} \beta_i \left(\mathcal{F}(\psi^{n - i}) - \psi^{n - i}\right) \right\|, &
    \text{subject to }\sum_{i = 0}^{m - 1} \beta_i &= 1,
\end{align}
and set $\psi^{n + 1} = \sum_{i = 0}^{m - 1} \beta_i \mathcal{F}(\psi^{n - i})$.

The iterations described in \Cref{alg: space-time} can be applied either per time step or over an entire simulation. Iterating between the flow and mechanics equations at each time step is often more efficient, particularly if the model is close to a steady state. However, it is also more invasive from an implementation perspective. 

On the other hand, iterating over an entire simulation is less invasive, but more memory demanding because the right hand sides of \eqref{eq: OPM solve} and \eqref{eq: TPSA solve} need to be saved and loaded. In the experiments of \Cref{sec:numerical_results}, we employ the latter iteration scheme.

Finally, if the dynamics of the system are sufficiently slow, then we may avoid iterating between the two systems. If we lag the influence of solid mechanics on the flow equations by one time step, then the problem effectively becomes a one-way coupled system. This leads us to the lagged scheme described by \Cref{alg: lagged}. 

\begin{algorithm}[H]
    \caption{Lagged scheme}
    \label{alg: lagged}
\begin{enumerate}[leftmargin=*,label=\arabic*.] 
    \item Initialize $p_f(t_0)$ and $\hat p(t_0)$ using the initial conditions, set $i = 0$, and $\hat p(t_{-1}) = \hat p(t_0)$. 
    \item At time $t = t_i$, use TPFA to solve the flow problem: find $p_f(t_{i + 1})$ such that
    \begin{align}
        (c_0 + \frac{\alpha^2}{\lambda}) \partial_t^{i + \frac12} p_f + \Mat{TPFA} p_f (t_{i + 1})
        &= f_p (t_{i + 1}) - \frac{\alpha}{\lambda} \partial_t^{i - \frac12} \hat p,
    \end{align}
    with $\partial_t^{i + \frac12}p \coloneqq \frac{p(t_{i + 1}) - p(t_i)}{\Delta t}$. 
    \item Use TPSA to solve for the solid mechanics variables
    \begin{align} 
        \Mat{TPSA} \begin{bmatrix}
                \bu(t_{i + 1}) \\ \br(t_{i + 1}) \\ \hat p(t_{i + 1})
            \end{bmatrix} =
            \begin{bmatrix}
                \bm{f_u}(t_{i + 1}) \\ 0 \\ - \frac{\alpha}{\lambda} \Delta p_f(t_{i + 1}))
            \end{bmatrix}
        \end{align}
    \item Increment $i$ and repeat the previous two steps until the end of the simulation.
\end{enumerate}
\end{algorithm}

\section{Implementation}
\label{sec:implementation}

We dedicate this section to the numerical implementation of \Cref{alg: space-time} and highlight how it can be used to easily introduce poroelasticity in existing flow simulators. \Cref{sub:solving_the_flow_equations} concerns the effects of the solid mechanics on the flow simulation and \Cref{sub:solving_the_linearized_elasticity_equations} proposes an efficient solver for the TPSA system.

\subsection{Solving the flow equations}
\label{sub:solving_the_flow_equations}

In order to implement the iterative coupling described in \Cref{alg: space-time} in existing numerical software, we need to consider two aspects.

First, we recognize the term $-\frac{\alpha}{\lambda}\partial_t \hat p$ on the right hand side of \eqref{eq: OPM solve} as a mass source. We incorporate this term in OPM Flow by using the keyword \texttt{SOURCE}, which allows for the prescription of an influx of fluid in each of the grid cells. 

Second, the term $\frac{\alpha^2}{\lambda} \partial_t p_f$ on the left hand side acts as an additional compressibility, cf.~\Cref{rem:compressibility}. We included this effect in OPM Flow by implementing a new keyword \texttt{ROCKBIOT}.


\subsection{Solving the TPSA system}
\label{sub:solving_the_linearized_elasticity_equations}

The TPSA system presented in Section \ref{sec:TPSA} is well-posed, but it may be challenging to solve numerically for two reasons. First, the scaling with material parameters is not favorable if the Lamé parameters are large. Secondly, we require efficient solvers in case the problem is too large for direct methods. In the following subsections, we propose propose left- and right-preconditioners to handle these challenges.

\subsubsection{Rescaling the system}

We first consider the dependency on material parameters. Let $Mx = b$ be short-hand notation for \eqref{eq:TPSA_system}
%
Inspecting the $3 \times 3$ block structure of $M$, we note that the $(1,1)$ block scales linearly with respect to $\mu$, the $(2,2)$ block scales as $\mu^{-1}$ and the $(3,3)$ block consists of two terms that scale as $\lambda^{-1}$ and $\mu^{-1}$, respectively. 

This mismatch in scaling may cause numerical problems if the parameters are large.
For example, granite rock has Lamé parameters on the order of $10 \unit{GPa}$ \cite{ji2010lame}. If the system is posed in Pascals, then the $\mu$ and $\mu^{-1}$ are in a ratio of $\sim10^{20}$, which leads to difficulties for the floating point arithmetic of the linear solver.

To counteract this imbalance between the equations, we propose a scaling of the rows and columns of the algebraic system. Let $\mu_0$ be the average of $\mu$ over $\Omega$, and let us assume that both $\frac{\mu_0}{\lambda}$ and $\frac{\mu_0}{\mu}$ are bounded from above by a reasonable constant.

We define the scaled matrix $\widetilde M$ and right-hand side $\tilde b$ by introducing the diagonal matrix $\Lambda \in \mathbb{R}^{n_{\text{dof}} \times n_{\text{dof}}}$ as follows:
\begin{align}
    \Lambda &= \begin{bmatrix}
        \mu_0^{-\frac12} \\
        & \mu_0^{\frac12} \\
        & & \mu_0^{\frac12}
    \end{bmatrix}, &
    \widetilde M &= \Lambda M \Lambda, &
    \tilde b &= \Lambda b.
\end{align}

Solving the system $Mx = b$ is equivalent to solving $\widetilde M\tilde x = \tilde b$ and retrieving $x = \Lambda \tilde x$. The system matrix $\widetilde M$ is easier to handle numerically because the scaling with $\mu$ has been removed from the diagonal blocks. The off-diagonal blocks, on the other hand, remain unchanged.

\subsubsection{Efficient preconditioning}
\label{subs:preconditioning}

The TPSA method employs a minimal stencil, which leads to a system matrix that is highly sparse. Nevertheless, for large-scale simulations, direct solvers are not feasible and we must rely on efficient and scalable iterative solvers. Since the system matrix of TPSA is not symmetric, we require a linear solver that handles general matrices. In our numerical experiments, the Bi-Conjugate Gradient Stabilized (BiCGStab) method \cite{van1992bi} proved to be well-suited.

As a Krylov subspace method, the performance of BiCGStab depends largely on the choice of preconditioner. We therefore propose a preconditioner in this subsection that is block-triangular and exploits the $3 \times 3$ block structure of the TPSA matrix. We consider the diagonal blocks in more detail.

First, the $(1, 1)$ block $M_{11}$ contains the finite volume discretization of $- \nabla \cdot (2 \mu \nabla \bu)$. This corresponds to three instances of a TPFA of a weighted Laplace problem, one for each component of $\bu$. To precondition such problems, we employ the Algebraic Multi-Grid (AMG) method, which is an effective preconditioner for Laplace problems discretized by TPFA. More specifically, we apply a single $V$-cycle of AMG by smoothed aggregation on the three independent sub-blocks of $M_{11}$ \cite{vanek1996algebraic}. We denote this operation by $\operatorname{AMG}_V(M_{11})^{-1}$.

Second, the $(2,2)$-block of the matrix $M$ is the discretization of $\mu^{-1} \br$. This is a diagonal matrix, so its inverse is directly available. 

The third and final block on the diagonal is a discretization of $\lambda^{-1}\hat p - \nabla \cdot (\delta^\mu \nabla \hat p)$ in which the second term forms the stabilization. It is thus composed of a diagonal matrix plus a two-point flux finite volume discretization of a Laplace-type operator. Again, we use a single $V$-cycle of AMG to approximate the inverse of this block, denoted by $\operatorname{AMG}_V(M_{33})^{-1}$.

Finally, we discard the upper-diagonal blocks of $M$, leading to the following, block-triangular preconditioner:
\begin{align} \label{eq: preconditioner}
    P \coloneqq 
    \begin{bmatrix}
        \operatorname{AMG}_V(M_{11}) \\
        M_{21} & M_{22} \\
        M_{31} & & \operatorname{AMG}_V(M_{33})
    \end{bmatrix}^{-1}
\end{align}

We emphasize that the matrix $P$ is not assembled. Instead, we realize the action of the linear operator by implementing the three-step forward substitution.

\section{Numerical results}
\label{sec:numerical_results}

In this section, we introduce three numerical test cases to examine various aspects of the coupling between TPSA and OPM. Convergence to an analytical solution is shown in \Cref{sub:spatial_convergence_to_an_analytical_solution} and we show the performance of the preconditioner from \Cref{subs:preconditioning} with respect to the mesh size. \Cref{sub:a_sealing_barrier} highlights the global effects that poroelasticity can introduce in a model that is otherwise compartmentalized and we compare the splitting schemes of \Cref{sub:iteration_schemes}. Finally, \Cref{sub:the_norne_geological_reservoir} investigates the performance of the coupled method on the complex geometry of the Norne geological reservoir.

The implementation details are as follows. The Python bindings for OPM Flow were implemented using \texttt{pybind11}, and the grid was handled using the package \texttt{opmcpg} \cite{opmcpg}. The AMG cycles in the preconditioner are implemented using \texttt{PyAMG} \cite{bell2022pyamg}. All numerical experiments were performed on a laptop with 16 processors (5\unit{GHz}) and 32\unit{GB} of RAM. The source code is available at \url{https://github.com/wmboon/tpysa}.

\subsection{Spatial convergence to an analytical solution}
\label{sub:spatial_convergence_to_an_analytical_solution}

Let $\Omega$ be a cube with side length of one \si{meter}. We define the fluid pressure and displacement by introducing an analytical function:
\begin{align} \label{eqs:analytical sol}
    \phi(\bx) &\coloneqq \prod_{i = 1}^3 \sin^2(\pi x_i), &
    p_f &\coloneqq \phi + p_0, &
    \bu &\coloneqq \sum_{i = 1}^3 \left(\partial_{x_{i + 1}} \phi - \partial_{x_{i - 1}} \phi \right) \bm{e}_i,
\end{align}
in which $p_0$ is the hydrostatic pressure due to gravity and 
$\bm{e}_i$ is the $i$-th canonical unit vector of $\mathbb{R}^3$. The indices in the subscripts are understood modulo 3. 

We moreover set the following material parameters. For the flow equations, we have a permeability of $K = 1$ \si{Darcy} and a viscosity of $\mu_w = 5 \times 10^{-4}$ \si{Pa.s}. The elastic medium has Lamé parameters $\mu = 0.01 $ \si{Pa} and $\lambda = 1$ \si{Pa}.
The boundary conditions are chosen to be no-flow and zero displacement. 
By substituting the solution in \eqref{eqs: Biot original}, we obtain analytical expressions for the body force $\bm{f_u}$ and mass source $f_p$.

We consider a sequence of Cartesian grids with decreasing mesh size $h$ and let the simulation reach a steady state by running 50 time steps with $\Delta t = 50$ \si{days}. The solution at final time is then compared to the analytical solution of \eqref{eqs:analytical sol}. 

\begin{figure}[ht]
    \centering
    \subfloat[]{
        \includegraphics[width=0.45\textwidth]{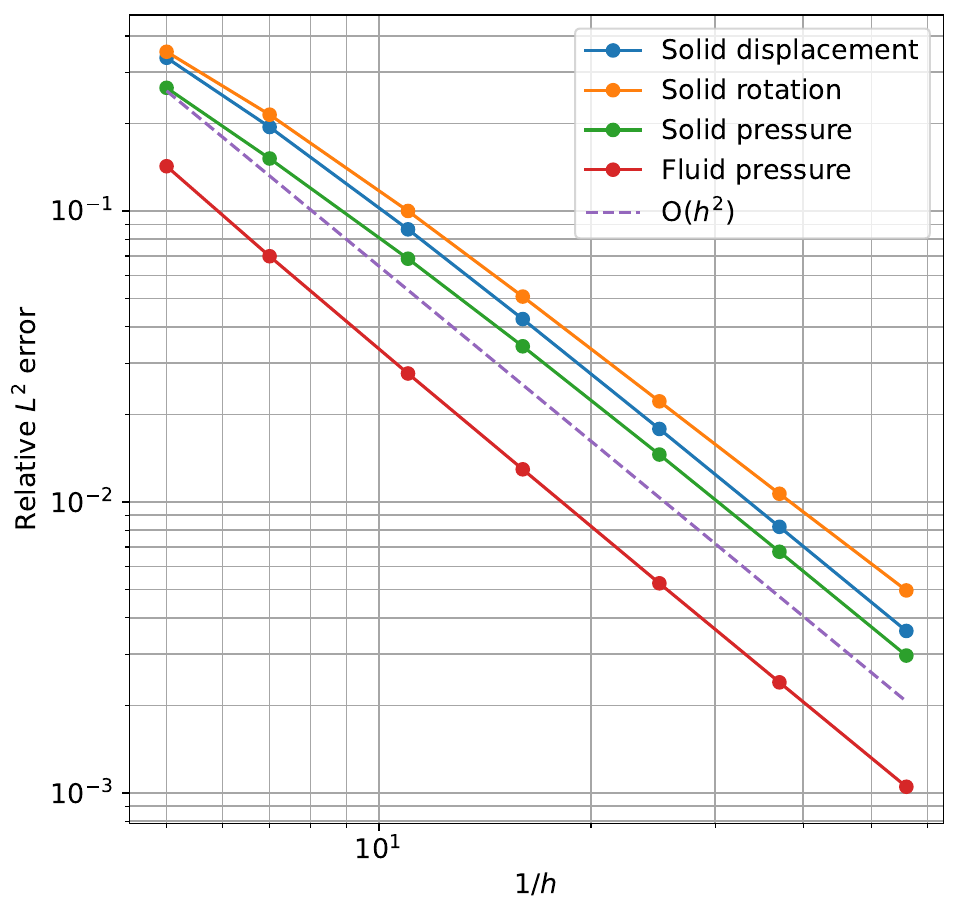} 
        \label{fig:L2_convergence}
        }
    \subfloat[]{
        \includegraphics[width=0.45\textwidth]{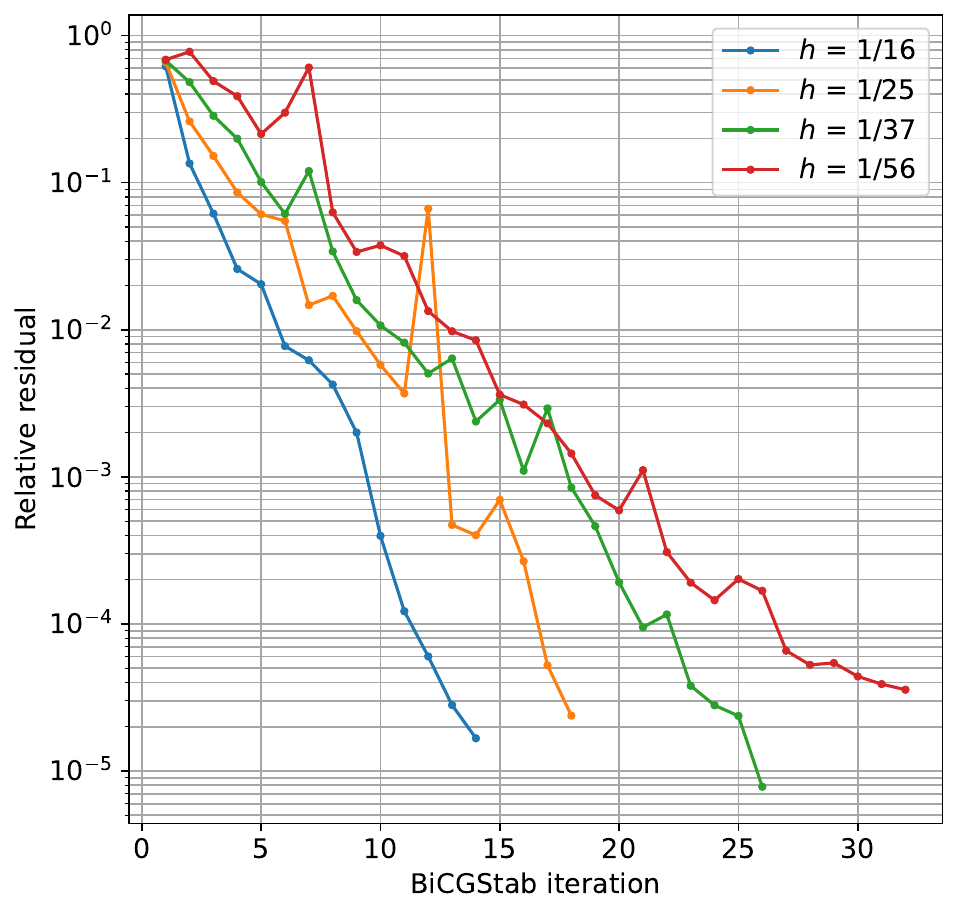} 
        \label{fig:BiCGSTab_analytical}
        }
    \caption{(a) The error converges quadratically with respect to the mesh size $h$ in all variables. (b) Convergence of the iterative solver BiCGStab for the TPSA system using the preconditioner \eqref{eq: preconditioner} on the four finest grids.}
\end{figure}

From the relative $L^2$ errors presented in \Cref{fig:L2_convergence}, we conclude that TPSA converges quadratically in all variables.  Theoretically, the method is guaranteed to converge only linearly for face-orthogonal grids and constant material parameters \cite{nordbotten2025two}. The increased convergence rate observed here is a consequence of the high regularity of the grids and the smoothness of the analytical solution.

We investigate the effectiveness of the preconditioner by considering the residual at each iteration of BiCGStab for the finest four grids, cf. \Cref{fig:BiCGSTab_analytical}. While the convergence is mainly monotone, we observe an increase in the required number of iterations as the mesh size decreases. 

In terms of computational cost, the iterative solver requires 2.4 \si{s} to solve the TPSA system with \num{109} \si{kdof} ($h = 1/25$), 6.1 \si{s} for \num{355} \si{kdof} ($h = 1/37$) and 24.5 \si{s} for $\num{1.23}$ \si{Mdof} ($h = 1/56$), on average per time step. These runtimes indicate a favorable scaling between the number of degrees of freedom and the solving time, which is mainly due to the parallelized software packages mentioned in the beginning of the section. Proper verification of the algorithmic scaling requires a more optimized implementation of the preconditioner, which is beyond the scope of this work.

\subsection{A sealing barrier}
\label{sub:a_sealing_barrier}

The second test case highlights a phenomenon that only occurs if poroelasticity is included in the model. In particular, we consider a compartmentalized system and show that poroelasticity introduces a global effect, even if the flow systems are mutually independent.

\begin{figure}[ht]
    \centering
    \subfloat[]{
        \includegraphics[width=0.45\textwidth]{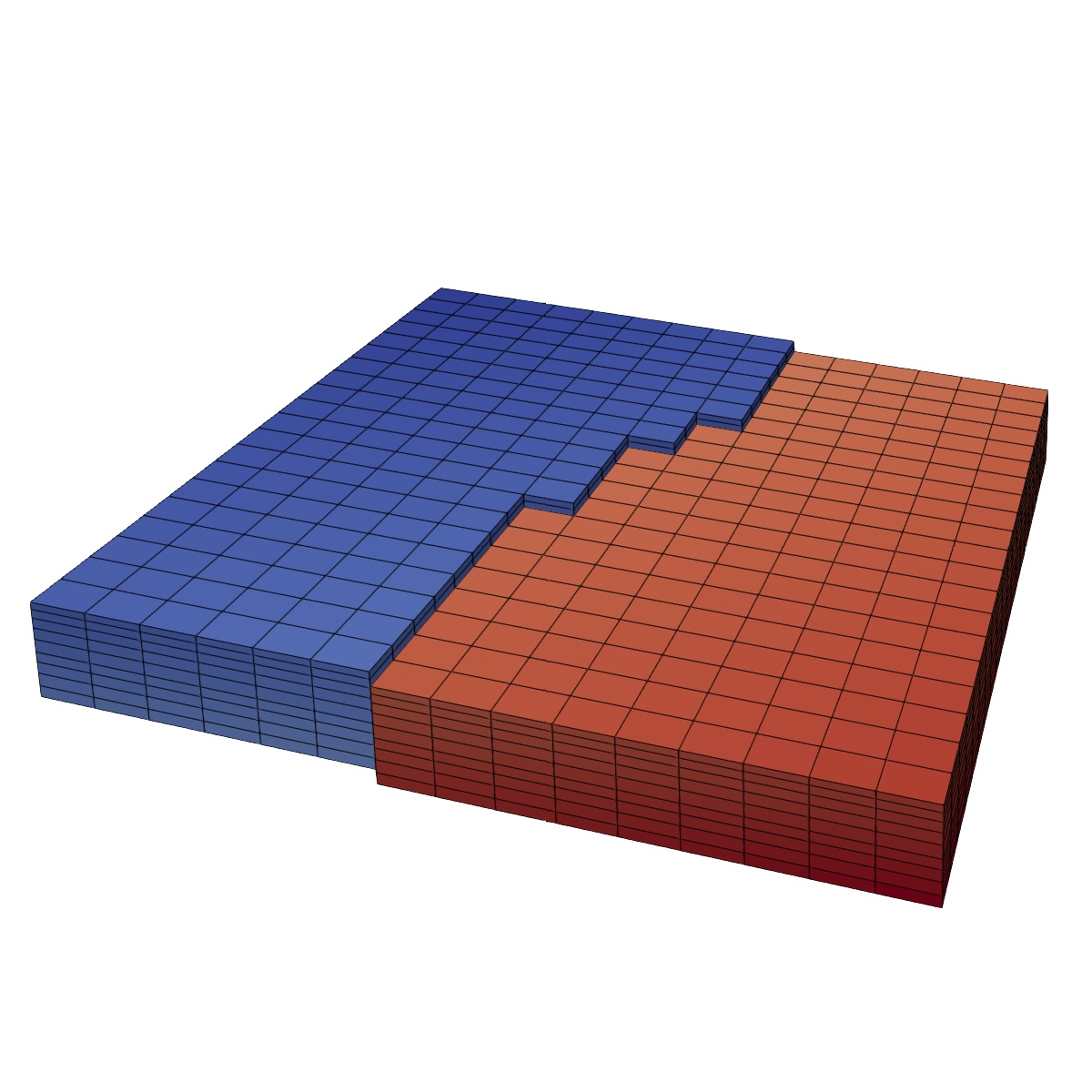}
        \label{fig:subdomains_fault}
    }
    \subfloat[]{
        \includegraphics[width=0.45\textwidth]{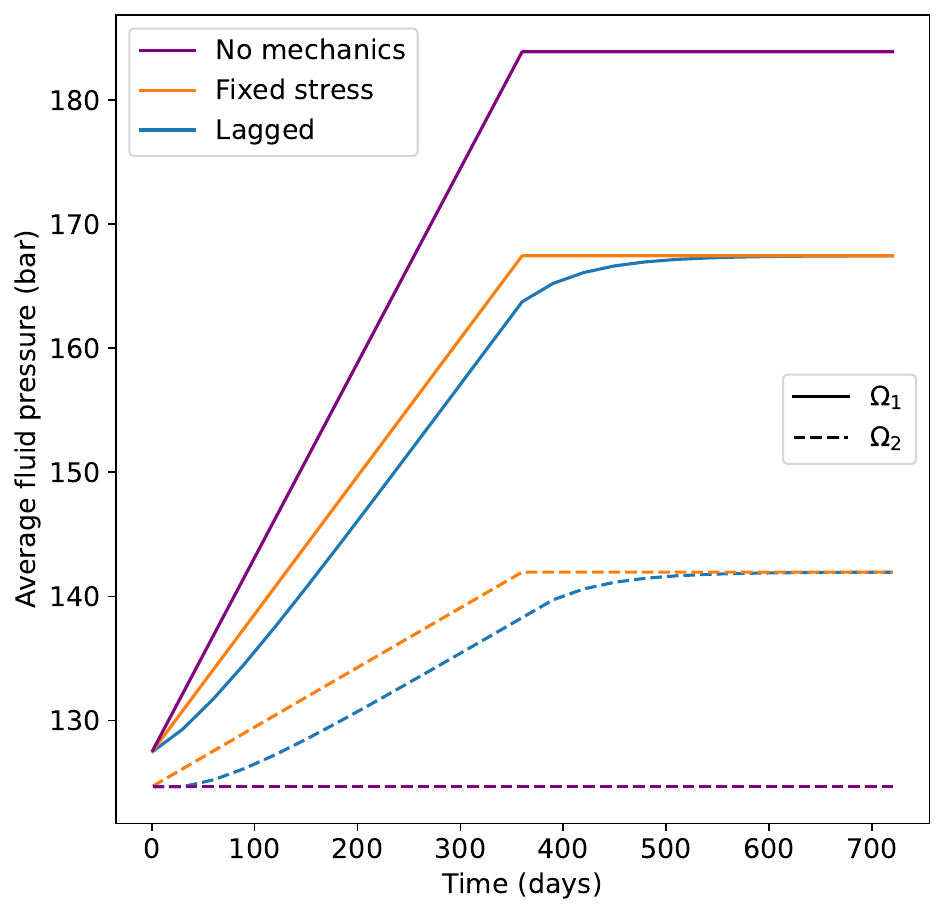}
        \label{fig:pressure_in_time}
    }
    \caption{(a) The second test case includes a sealing barrier that divides the domain into two subdomains, $\Omega_1$ in red and $\Omega_2$ in blue, respectively. (b) The fluid pressure, averaged over the subdomains. The pressure increases due to an injection well in $\Omega_1$, which only affects the pressure in $\Omega_2$ if poroelasticity effects are included in the model.}
\end{figure}

Let the domain $\Omega$ be subdivided into two subdomains by an impermeable barrier as illustrated in \Cref{fig:subdomains_fault}. We impose no-flow and zero displacement conditions on all boundaries. An injection well is introduced in subdomain $\Omega_1$ that injects a compressible fluid with a rate of $100\unit{m^3/day}$ for the first 360 days of the simulation (12 time steps). The injection is then stopped, letting the reservoir equilibrate for another 360 days.

If solid mechanics effects are neglected, then this problem decouples into two independent flow problems on the respective subdomains. The well increases the pressure in $\Omega_1$ while the pressure in $\Omega_2$ remains unaffected, as reflected by the purple curves in \Cref{fig:pressure_in_time}.

However, if poroelastic effects are incorporated, then the medium is allowed to deform. The increased pressure near the injection well causes subdomain $\Omega_1$ to expand. In turn, the barrier bulges since the boundaries are clamped, decreasing the volume of $\Omega_2$ which, in turn, causes the fluid pressure to increase there. This effect is clearly visible in the dashed pressure curves of $\Omega_2$ in \Cref{fig:pressure_in_time}. We remark that a simpler model, in which the rock compressibility is included in the storativity coefficient $c_0$, cannot capture this inter-subdomain pressure influence.

Moreover, the inclusion of poroelasticity causes the pressure in $\Omega_1$ to reach a lower steady state. The fluid pressure is thus effectively dissipated into a region that is not connected by flow paths. We can make this observation more exact by integrating the mass balance equation \eqref{eq: mass balance orig} and applying the boundary and initial conditions:
\begin{align} \label{eq: mean pressure}
    \int_{0}^T\!\!\int_\Omega f_p \ \myd \bx \myd t 
    &= \int_{0}^T\!\!\int_\Omega \partial_t (c_0 p_f + \alpha \nabla \cdot \bu) + \nabla \cdot \bq \ \myd \bx \myd t \nonumber\\
    &= 
    c_0 \int_\Omega \Delta p_f(T) \ \myd \bx + \int_0^T\!\!\int_{\partial \Omega} \partial_t \alpha \bn \cdot \bu + \bn \cdot \bq \ \myd \bm{s} \myd t \nonumber\\
    &= 
    c_0 \int_\Omega \Delta p_f(T) \ \myd \bx
\end{align}

Equation \eqref{eq: mean pressure} shows that the average pressure deviation only depends on the cumulative mass added to the system. Hence, by allowing the pressure to increase in $\Omega_2$ through mechanical effects, the pressure in $\Omega_1$ reaches a lower equilibrium state.

The phenomenon extends from sealing barriers, showcased here, to impermeable regions in the model. Such regions can similarly connect different areas of the reservoir through the solid mechanics equations, and thereby play a significant role in the fluid pressure evolution. This indicates that care must be taken when removing regions with low permeabilities or low pore volumes from the reservoir simulation, even if they do not facilitate any significant fluid flow.

Secondly, \Cref{fig:pressure_in_time} illustrates how the choice of splitting schemes affects the pressure evolution in time. The Lagged scheme (\Cref{alg: lagged}) clearly exhibits a delay in the pressure response compared to the converged fixed stress solution. The steady-state solution at the end of the simulation is the same, however, and we have only observed significant differences between the schemes if large shocks are induced, e.g. through the boundary conditions. 

\begin{figure}[ht]
    \centering
    \includegraphics[width=0.45\textwidth]{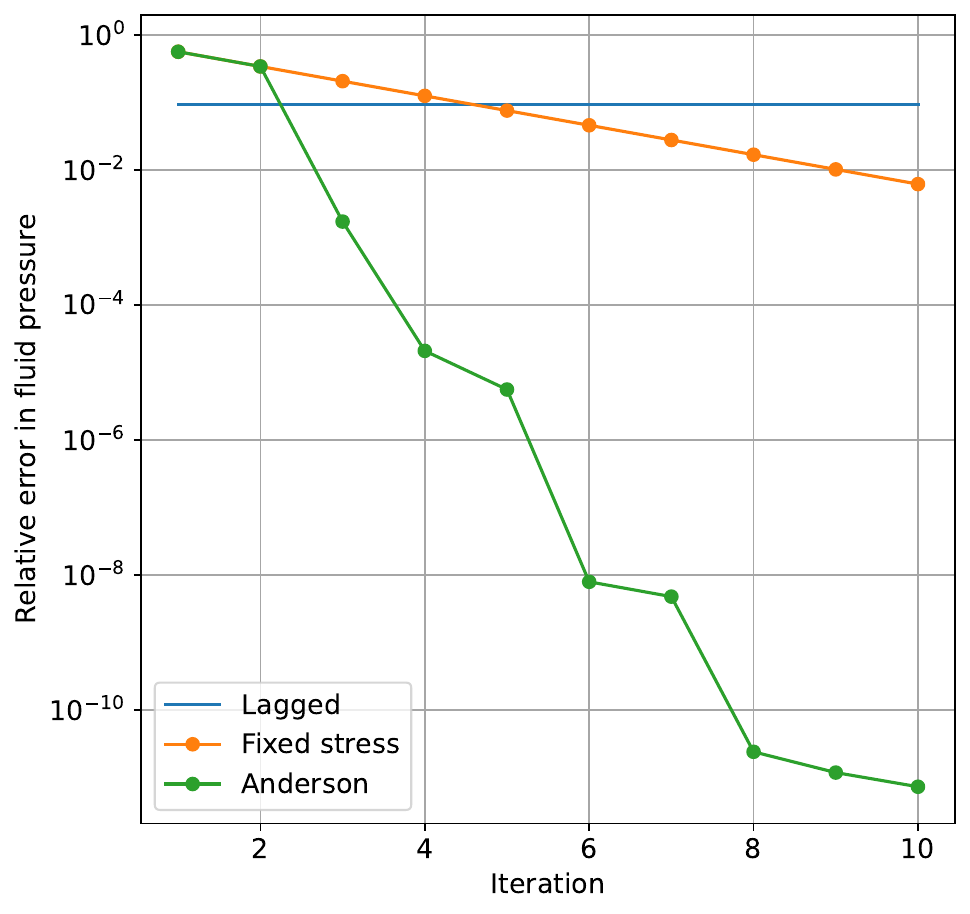}
    \caption{Comparison between the different splitting schemes for the test case of \Cref{sub:a_sealing_barrier}. We observe that Anderson acceleration is particularly effective for this simple problem.}
    \label{fig:fixedpoint_vs_lagged_fault}
\end{figure}

The convergence of the splitting schemes from \Cref{sub:iteration_schemes} is illustrated in \Cref{fig:fixedpoint_vs_lagged_fault}. \Cref{alg: space-time} is realized as follows. At iteration $n$, the mass source $\psi^n = - \frac{\alpha}{\lambda} \partial_t \hat p^{n - 1}$ is provided for all $t \in [0, T]$, and the resulting $\psi^{n + 1}$ becomes the input for the next iteration. Thus, each iteration requires a complete simulation of the problem. 

The introduction of Anderson acceleration significantly improves the convergence of the splitting scheme. The first iteration is the same as in the fixed point iteration because the scheme starts with a zero initial guess for $\psi$. The second iteration also coincides because the solution to \eqref{eq: Anderson} is trivially $\bm{\beta} = \beta_0 = 1$. From the third iteration onwards, the accelerated scheme immediately outperforms the simple fixed point scheme. The additional cost consists of solving the $m \times m$ least-squares problem \eqref{eq: Anderson}, which is negligible for our choice of $m \le m_0 \coloneqq 5$.

Finally, we remark that the grid is in a corner-point format and the two subdomain grids are non-matching across the barrier. Faces at the barrier are subdivided according to the overlaps between cell boundaries. In turn, adjacent cells have more than six faces, i.e. each becomes a polyhedron that is more general than a hexahedron. This processing of the grid is already performed for the TPFA flow discretization, and we reuse this geometric information in the TPSA assembly.

The grid consists of \num{2700} cells, which leads to \num{18900} degrees of freedom for TPSA. The system matrix is thus sufficiently small to save the LU-factorization of the system matrix (3 sec), yielding fast mechanics solves at each time step. For completeness, we also tested BiCGStab with preconditioner \eqref{eq: preconditioner}, which reaches a relative residual of $10^{-5}$ within 25 iterations (0.3 sec), at each time step. 

\subsection{The Norne geological reservoir}
\label{sub:the_norne_geological_reservoir}

In the third test case, we consider a challenging geometry to showcase the flexibility of the coupled finite volume method. We choose the Norne reservoir as our domain of computation, for which a conforming corner-point grid is available. Since no analytical solution is available, nor a sequence of grids, we mainly present qualitative observations in this section.

We impose no-flux boundary conditions for the flow system. Let the Biot-Willis constant be $\alpha = 0.87$ and the Lamé parameters $\mu = 3.5$ \si{GPa} and $\lambda = 4.0$ \si{GPa}, in accordance with \cite{ji2010lame}. For the solid mechanics, we impose Robin, or spring, boundary conditions with $w_j \coloneqq \frac{0.05 \Delta z}{\mu}$ in which $\Delta z \approx 620$ \si{m} is the maximal vertical extent of the geometry. This effectively surrounds the domain by springs that are each $0.05 \Delta z$ long and have the same stiffness as the reservoir. In turn, the domain boundary is able to deform, with larger deformations inducing larger opposing forces. A discussion regarding the implementation of Robin boundary condition is given at the end of \Cref{sub:TPSA_discretization}.

At the start of the simulation, a water injection well starts operating at a rate of $10^3$ \si{m^3/day}. Simultaneously, an extraction well starts in the rear of the domain with the same rate. We let the wells operate for 15 time steps of $\Delta t = 20$ \si{days}. The wells are then shut and we let the system equilibrate for a second set of 15 time steps.

\begin{figure}[ht]
    \centering
    \subfloat[]{
        \includegraphics[width=0.45\textwidth]{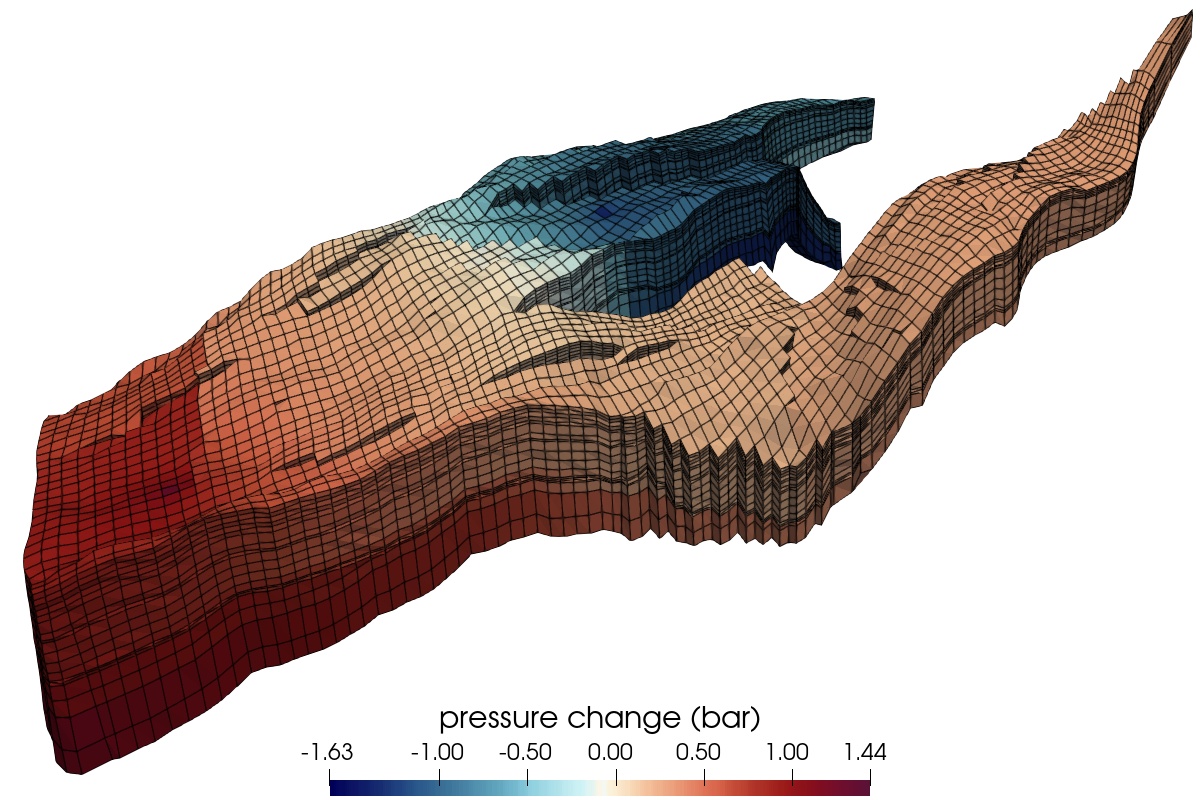}
        \label{fig:pressure_Norne}
    }
    \subfloat[]{
        \includegraphics[width=0.45\textwidth]{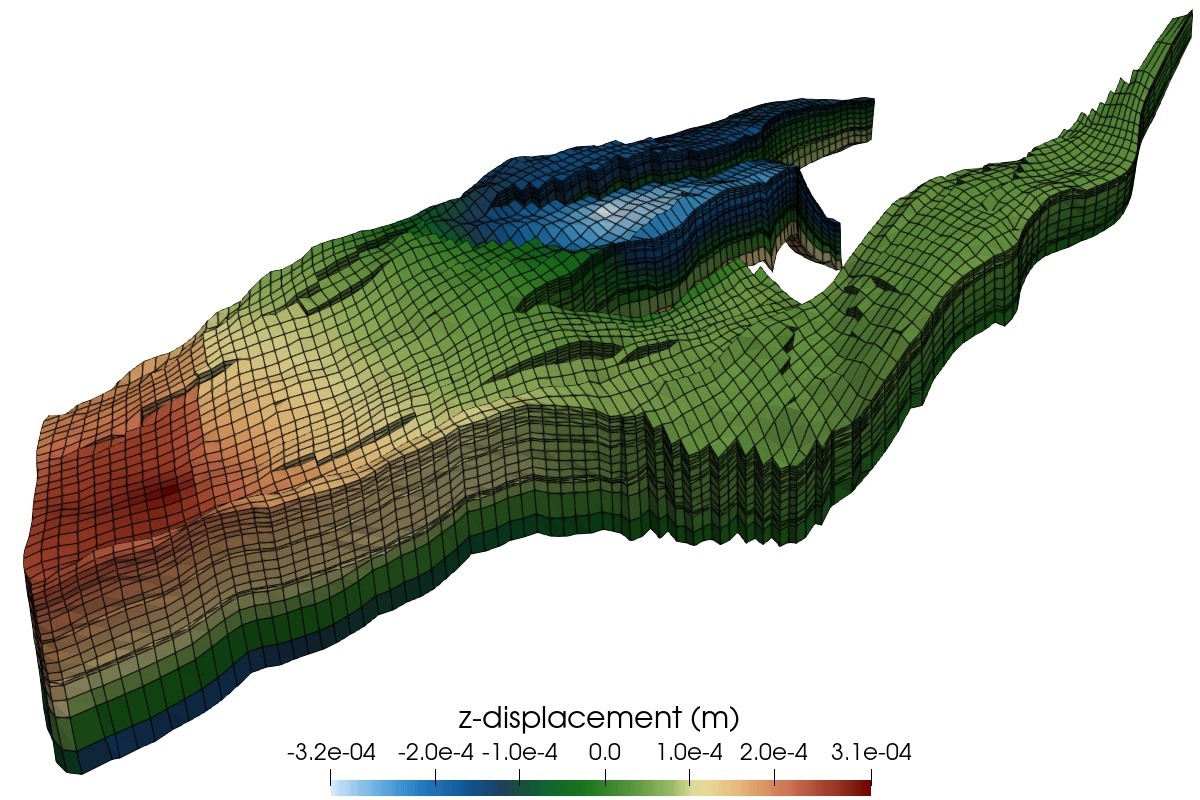}
        \label{fig:displacement_Norne}
    }
    \caption{(a) Deviation in fluid pressure $\Delta p_f$ and (b) the vertical displacement $u_z$ for the third test case on the Norne reservoir. At the injection well in the lower left of the region, the fluid pressure increases, causing the porous medium to expand. Conversely, a decrease in pressure makes the reservoir contract at the extraction well in the rear.}
    \label{fig:Norne}
\end{figure}

\Cref{fig:Norne} presents the solution when the wells are shut down, at $t = 300$ \si{days}. Because of the high aspect ratio of the reservoir, we have amplified the domain in the $z$-direction by a factor 5. 

In \Cref{fig:pressure_Norne}, we clearly see a higher pressure in the lower left of the domain, where the injection well is located. The extraction well, on the other hand, causes a pressure drop in the rear regions of the model. While the pressure is fairly continuous, we note several discontinuities, or jumps, which are a consequence of lower transmissibilities between geological regions defined in the model. 

The resulting displacement field is presented in \Cref{fig:displacement_Norne}. We clearly observe that the medium expands near the injection well. In particular, the top of the domain is lifted and the bottom is pushed downward. Conversely, the reservoir contracts near the extraction well, visible by the reversed color gradient in the vertical direction. 

\begin{figure}[ht]
    \centering
    \subfloat[]{
        \includegraphics[width=0.45\textwidth]{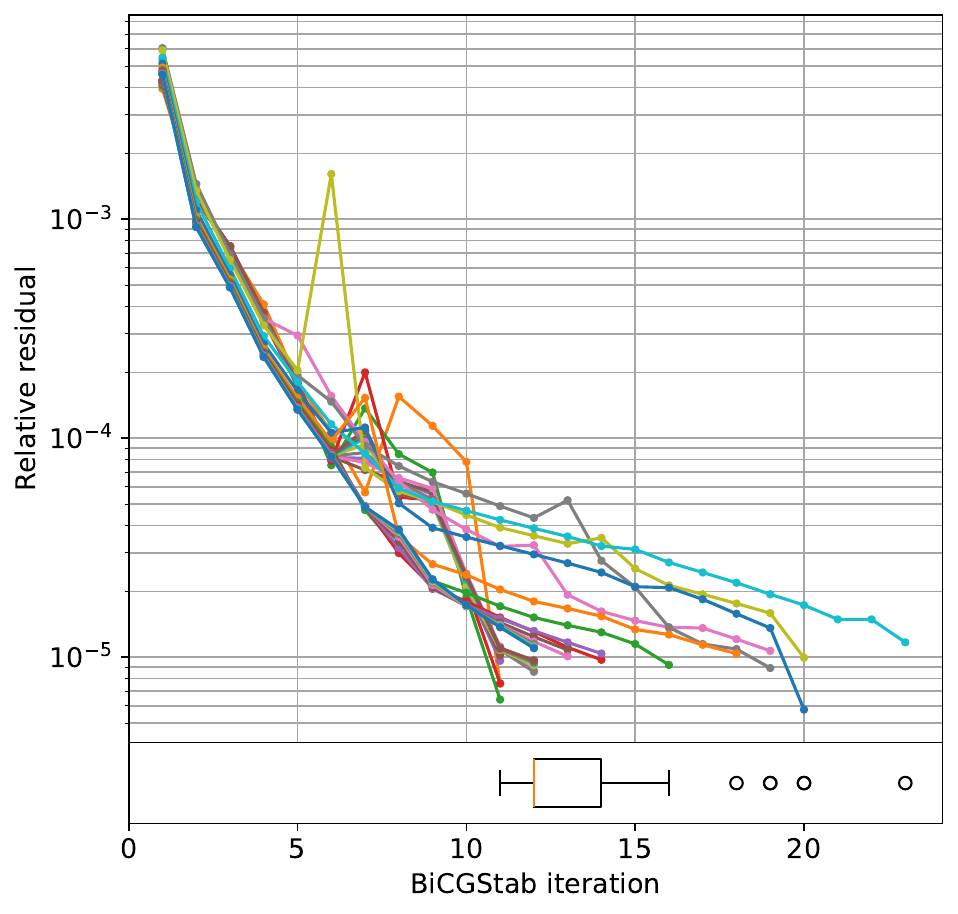}
        \label{fig:BiCGStab_Norne}
    }
    \subfloat[]{
        \includegraphics[width=0.45\textwidth]{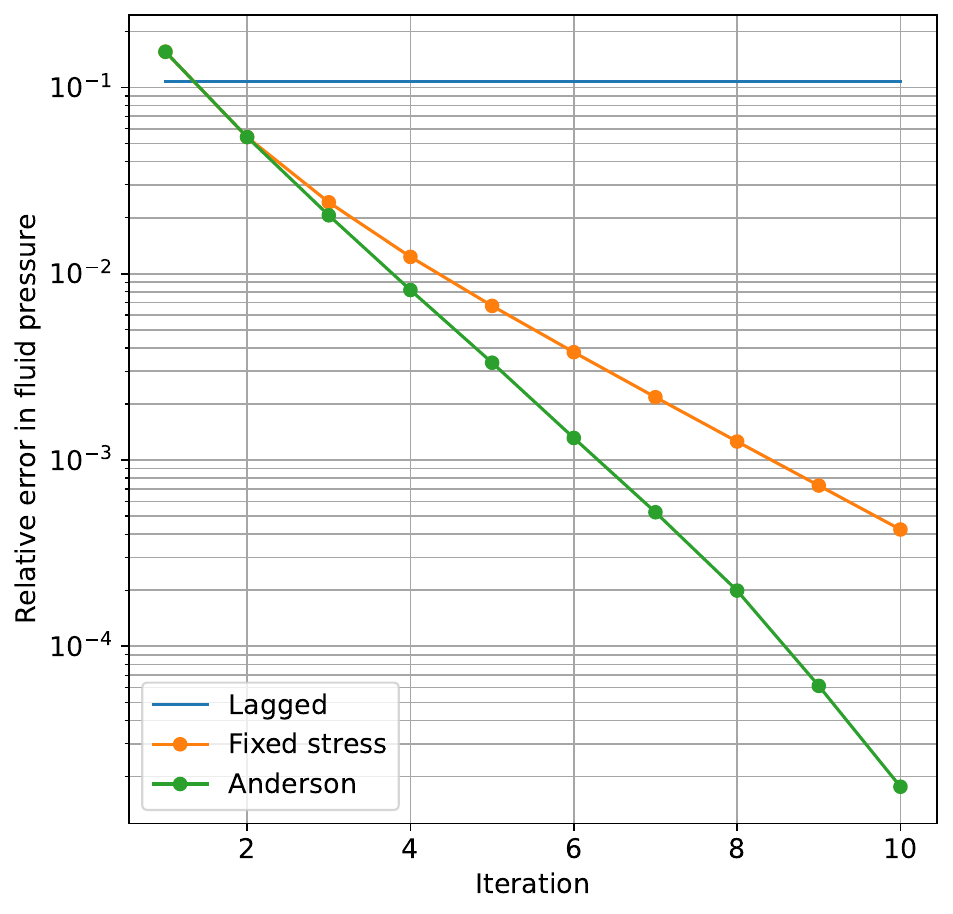}
        \label{fig:fixedpoint_vs_lagged_Norne}
    }
    \caption{(a) The proposed preconditioner \eqref{eq: preconditioner} allows BiCGStab to converge within 16 iterations for the majority of time steps. (b) The fixed stress schemes outperform the Lagged coupling from the second iteration onwards. The convergence improves if Anderson acceleration is applied. }
\end{figure}

Next, we report on the computational cost. The grid consists of \num{38180} cells, resulting in approximately 267kdof for TPSA. The TPSA system \eqref{eq:TPSA_system} was assembled and the AMG preconditioner \eqref{eq: preconditioner} was initialized in approximately one second. The system and preconditioner were saved and used at each time step.

\Cref{fig:BiCGStab_Norne} shows the performance of the preconditioner through the iterative solver BiCGStab, which reached the desired tolerance level within 16 iterations for the majority of time steps, often terminating after only 12. The mechanics solution was obtained in approximately two seconds per time step, on average. 

Finally, we compare the different poroelasticity coupling schemes in \Cref{fig:fixedpoint_vs_lagged_Norne}. The lagged scheme does not involve any iterations, and it therefore serves as an efficient baseline alternative. The simple fixed point iterations of \Cref{alg: space-time} already produces a more accurate pressure distribution on the second iteration. We emphasize that this corresponds to two full simulations of the problem. Again, we start seeing the positive effects of Anderson acceleration at the third iteration. After ten iterations, the error produced by the accelerated scheme is an order of magnitude smaller than the straightforward fixed point approach.

\section{Conclusions}
\label{sec:conclusions}

We proposed a solver for Biot poroelasticity that combines industrial finite volume solver for flow with the recently developed two-point stress approximation method for linearized elasticity. 
Through a reformulation of the problem, we propose a non-invasive fixed stress coupling scheme that utilizes tailored solvers for the two subproblems. Due to the coupling of two finite volume codes, the information passed between the solvers is minimal, as it concerns one value per cell, per time step. 

Using numerical experiments, we showcased the performance and flexibility of the approach. Spatial convergence was verified and we showed that the proposed preconditioner scales favorably with respect to the number of degrees of freedom. The fixed stress coupling scheme converges monotonically and its performance improved significantly after applying Anderson acceleration.

In the second test case, we highlighted the need for including poroelastic effects in a reservoir simulation. In particular, we showcased how the elasticity equations affect the fluid pressure evolution in otherwise isolated regions of the domain. 
This experiment shows that caution should be exercised before disregarding regions in the domain with low permeabilities or small pore volumes. While these regions may not facilitate significant fluid flow, they can still form important mechanical connections within the domain. Particularly for applications in which the maximum pressure is important, this additional dissipation of the pressure may have a significant impact.

\backmatter
\bmhead{Acknowledgments}

The MuPSI project (Multiscale Pressure-Stress Impacts on fault integrity for multi-site regional CO2 storage) is awarded through the Clean Energy Transition Partnership (CET-P) project number CETP-FP-2023-00298, with funding provided by the RCN Research Council of Norway, Scottish Enterprise, NWO Dutch Research Council, AEI-Agencia Estatal de Investigación, and US DoE, with contributions from Storegga Ltd,  Equinor ASA, Norske Shell AS, and EBN Capital BV.

\bmhead{Declarations}

The authors have no relevant financial or non-financial interests to disclose.

\bibliography{references}

\end{document}